\documentclass[12pt]{article}
\usepackage{geometry}
\usepackage{amsmath}
\usepackage{amsthm}
\usepackage{amssymb}
\usepackage{graphicx}
\graphicspath{ {./image/} }
\usepackage{enumitem}
\usepackage{caption}
\usepackage{subcaption}

\usepackage{tocloft}
\usepackage[hidelinks]{hyperref}

\usepackage{mathtools, nccmath}

\usepackage{cite}
\makeatletter
\newcommand{\citecomment}[2][]{\citen{#2}#1\citevar}
\newcommand{\citeone}[1]{\citecomment{#1}}
\newcommand{\citetwo}[2][]{\citecomment[,~#1]{#2}}
\newcommand{\citevar}{\@ifnextchar\bgroup{;~\citeone}{\@ifnextchar[{;~\citetwo}{]}}}
\newcommand{\citefirst}{\@ifnextchar\bgroup{\citeone}{\@ifnextchar[{\citetwo}{]}}}

\makeatother

\newtheorem{theorem}{Theorem}[section]

\newtheorem{conj}[theorem]{Conjecture}
\newtheorem{definition}[theorem]{Definition}
\newtheorem{?}[theorem]{Question}
\newtheorem{remark}{Remark}[section]

\begin{document}
	
\title{Bipyramid Volume, Mahler Measure and Some $\mathbb{Z}^2$-periodic Links}

\author{Hong-Chuan Gan}

\date{}

\maketitle              

\begin{abstract}
Champanerkar, Kofman and Lal\'{i}n conjectured an inequality between bipyramid volume of links and Mahler measure of associated dimer models induced from alternating links on torus.
Hyperbolic volume and Mahler measure can be related for isoradial graphs, which allows us to confirm the conjecture for two examples.
By exploiting a connection between perfect matchings of dimer models and spanning trees on lattices, five more examples are calculated.
\end{abstract}

\tableofcontents

\section{Introduction}
A knot or link can be represented as a diagram in the plane.
Graphs can be constructed from this diagram, for example the projection graph and the Tait graphs.
One can consider statistical models on these graphs and find interesting connections between quantities from knots and that from statistical models.
See \cite{jones1990knot} for a nice introduction of Jones to this topic.

The Vol-Det conjecture of Champanerkar-Kofman-Purcell\cite[Conjecture 1.10]{champanerkar2016geometrically} and the Champanerkar-Kofman-Lal\'{i}n conjecture\cite[Conjecture 1]{Champanerkar_2018} fall into this category---the interplay between knot theory and statistical mechanics.
The Vol-Det conjecture suggests an inequality between the hyperbolic volume of a hyperbolic alternating knot $K$ and its determinant $\det(K)$
$$ vol(K) < 2\pi \log \det (K) .$$
Here $\det(K) =  |\Delta(-1)| = |V(-1)|$ where $\Delta(t)$ is the Alexander polynomial and $V(t)$ is the Jones polynomial of $K$.
The Champanerkar-Kofman-Lal\'{i}n conjecture(called CKL conjecture below) is closely related to the Vol-Det conjecture, see Section \ref{section::ckl_construction} for detailed explanations.
It suggests an inequality between the bipyramid volume of an alternating link on torus and the Mahler measure of the characteristic polynomial of the dimer model induced from it.
See Section \ref{section::CKL conjecture} for definitions.
The CKL conjecture is the main motivation for this work.
We refer readers to \cite{Champanerkar_2018} and references therein for more motivations of  the conjecture.
Their examples are very helpful for this work.

Here we verify the CKL conjecture for some examples.
In Section \ref{section::CKL conjecture} we confirm the CKL conjecture for two examples using a connection between Mahler measure and hyperbolic volume for isoradial graphs explained in Section \ref{section::Hyperbolic volume and Mahler measure}.
In Section \ref{section::spanning tree entropy}, we use results on spanning tree entropy to calculate 5 more examples.
\subsection*{Acknowledgements.}
I would like to thank Abhijit Champanerkar and Ilya Kofman for correcting an error in the computation.
See Remark \ref{remark::error} for more details.

\section{Background on dimers} \label{section::background}
See \cite{kenyon2006dimers} for more details.
Let $G$ be a $ \mathbb{Z}^2 $-periodic bipartite planar graph.
Let $G_n$ be the quotient of $G$ by the action of $n\mathbb{Z}^2$.
It is a finite bipartite graph on a torus.
A \textit{dimer cover} or \textit{perfect matching} of a graph $G$ is a set of edges with the property that every vertex is contained in exactly one edge of the set.
If $\nu: E(G) \to \mathbb{R}_{>0}$ is a positive weight function on edges of $G$, we associate a weight $\nu(m) = \prod_{e\in m} \nu(e)$ to a dimer cover which is the product of its edge weights.
Define the \textit{partition function} of $G_1$ to be $Z=\sum_{m\in M} \nu(m)$ where $M$ is the set of all dimer covers of $G$.
\begin{description}[leftmargin=0pt]
	\item[Kasteleyn matrix.] 
	A Kasteleyn matrix for a finite bipartite planar graph $\Gamma$ is a signed, weighted adjacency matrix, with rows indexed by the white vertices and columns indexed by the black vertices, with $K(w, b)=0$ if $w$ and $b$ are not adjacent, and $K(w, b)=\pm \nu(wb)$ otherwise, where the signs are chosen so that the product of signs around a face is $(-1)^{k+1}$ for a face of degree $2k$. 
	Kasteleyn \cite{kasteleyn1967graph} showed that $|\det K|$ is the partition function $|\det K|  = Z(\Gamma)$.
	
	\item[Characteristic polynomial.] 
	For a finite bipartite graph on a torus, we construct a ``magnetically altered'' Kasteleyn matrix $K(z, w)$ as follows.
	Let $\gamma_x$ (resp. $\gamma_y$) be a path in the dual of $G_1$ winding once horizontally (resp. vertically) around the torus. 
	Multiply each edge crossed by $\gamma_x$ by $z^{\pm 1}$ depending on whether the black vertex is on the left or on the right, and similarly for $\gamma_y$.
	Define the \textit{characteristic polynomial} of $G$ by $P(z, w) = \det K(z, w)$.
	The \textit{Ronkin function} of $P(z, w)$ is
	$$ F(X, Y) = \frac{1}{(2\pi i)^2} \int_{\mathbb{T}^2} \log | P(e^X z, e^Y w) | \frac{dz}{z}\frac{dw}{w}. $$
	Note that $ F(0,0) $ is the \textit{Mahler measure} of $P(z, w)$(denoted as $M(P(z,w))$).

	The partition functions for graphs on a torus can be expressed in terms of the characteristic polynomial:
	$$ Z=\frac{1}{2}(-P(1,1)+P(1,-1)+P(-1,1)+P(-1,-1)). $$
	
	\item[Gauge equivalence.] 
	We define gauge equivalence following \cite{kenyon2018dimersandcirclepatterns}.
	Two weight functions $\nu_1, \nu_2$ are said to be \textit{gauge equivalent} if there is a function $F: V(\mathcal{G}) \to \mathbb{R}$ such that for any edge $vv'$, $ \nu_1(vv') = F(v)F(v')\nu_2 (vv') $.
	For a planar bipartite graph, two weight functions are gauge equivalent if and only if their \textit{face weights} are equal, where the
	face weight of a face with vertices $w_1, b_1, . . . , w_k, b_k$ is the ``alternating product" of its edge weights,
	$$ X = \frac{\nu(w_1b_1)\cdot \cdot \cdot \nu(w_kb_k)}{\nu(b_1w_2) \cdot \cdot \cdot \nu(b_kw_1)}. $$
	Different choices of gauge correspond to multiplying $K$ on the right and left by diagonal matrices with positive diagonal entries.
	Thus scaling edge weights of all edges incident to a vertex by a factor of $s$ would change $\log Z$ to $\log Z + \log s$.
\end{description}

\section{Hyperbolic volume and Mahler measure}\label{section::Hyperbolic volume and Mahler measure}
In this section, we describe how hyperbolic volume and Mahler measure can be related for isoradial graphs.
\subsection{Kenyon's formula} \label{section::kenyon's formula}
\begin{figure}
	\begin{subfigure}{0.5\textwidth}
		\centering
		\includegraphics[width=0.8\linewidth]{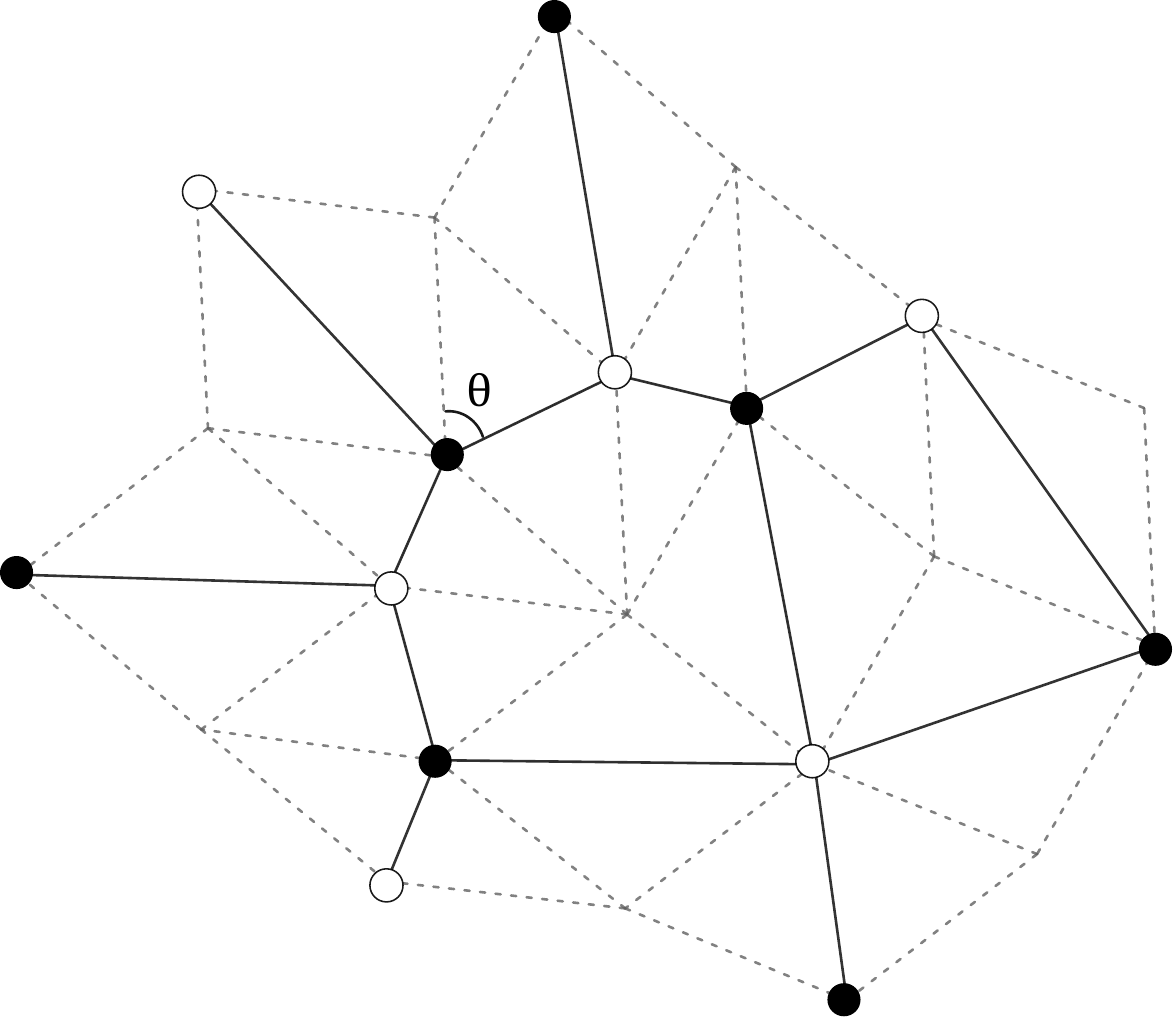}
		\caption{}
		\label{img::isoradialGraphLeft}
	\end{subfigure}%
	\begin{subfigure}{0.5\textwidth}
		\centering
		\includegraphics[width=0.8\linewidth]{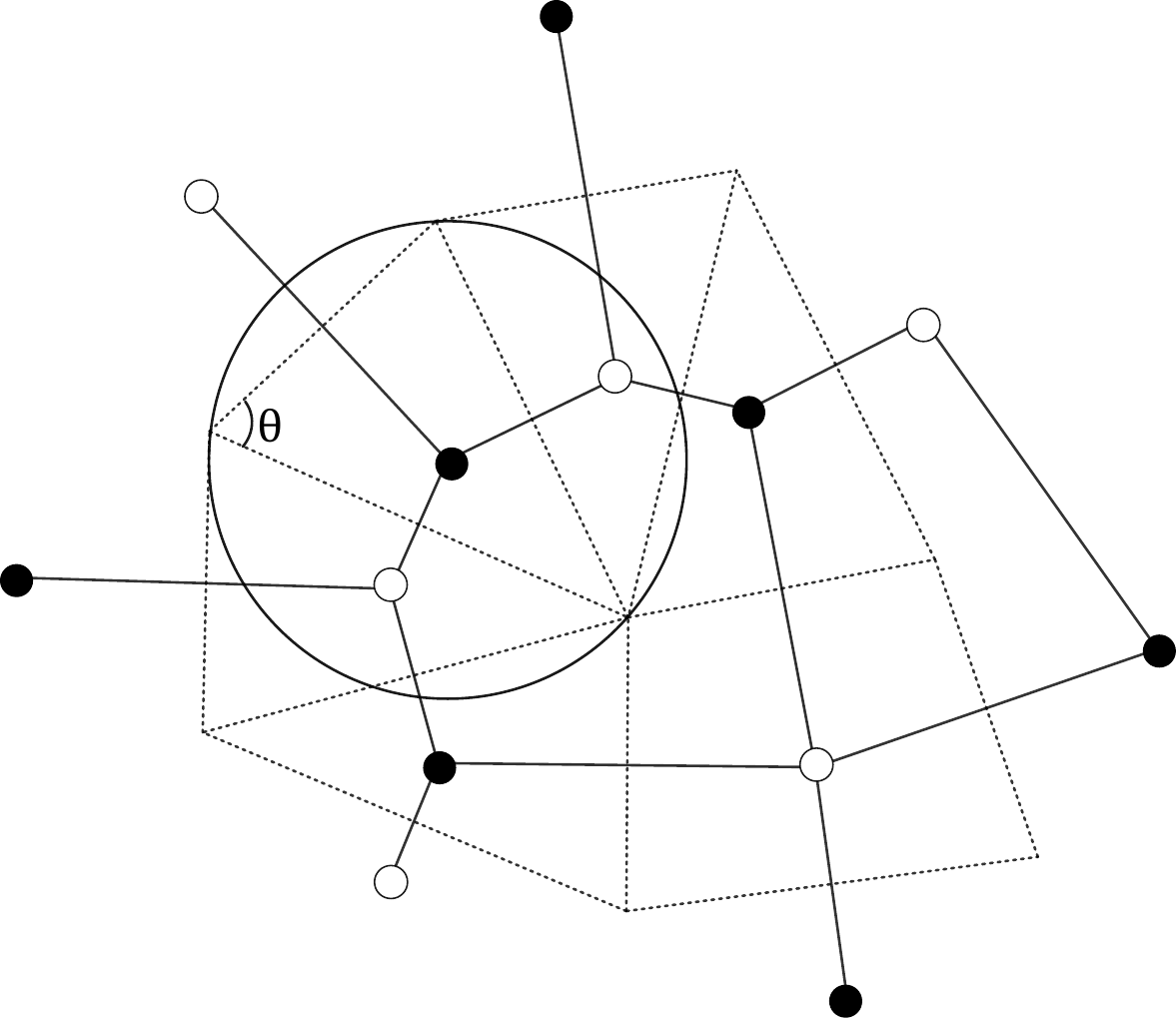}
		\caption{}
		\label{img::isoradialGraphRight}
	\end{subfigure}
	\caption{Left: an isoradial bipartite graph and the rhombus tiling. Right: an isoradial bipartite graph and the dual graph.}
\end{figure}

\begin{definition}
	An isoradial embedding of a planar graph is an embedding satisfying that
	\begin{itemize}
		\itemsep0em
		\item every face is a cyclic face,
		\item circumcenters of the faces are contained in the closures of the faces,
		\item all circumcircles have the same radius.
	\end{itemize} 
\end{definition}
Given the isoradial embedding $\mathcal{G}_I$ of $\mathcal{G}$, one can construct a graph called \textit{rhombus tiling} where vertices are given by circumcenters of circumcircles of faces of $\mathcal{G}_I$ and edges are given by connecting the circumcenter and vertices of the cyclic face for every face of $\mathcal{G}_I$.
Rhombus tiling gives a rhombus for every edge of $\mathcal{G}$.
We associate the number $\nu (e) = 2 \, \sin \, \theta_e$ to the edge $e$ where $2\theta_e$ is the angle of the coresponding rhombus in the rhombus tiling at the vertices it has in common with $e$.
This gives a function $\nu : E(\mathcal{G}_I) \to \mathbb{R}$ called \textit{critical weight function}.
See Figure \ref{img::isoradialGraphLeft}.
Note that $\theta_e$ is also the angle cut by the dual edge in the circumscribed circle of the dual graph.
See Figure \ref{img::isoradialGraphRight}.

Kenyon introduced the Dirac operator $\overline{\partial}$ and its normalized determinant $ {\det}_1 \overline{\partial} $(we refer readers to \cite{kenyon2002laplacian} for the definitions) and gave the following formula\cite[Theorem 1.2]{kenyon2002laplacian} relating $ {\det}_1 \overline{\partial}$ and the hyperbolic volume for isoradial embedding of a periodic bipartite graph with critical edge weights $\nu$
$$
\log \; {\det}_1 \overline{\partial} = \frac{2}{N} \sum_{edge \; e} \frac{1}{\pi}L(\theta_e)+\frac{\theta_e}{\pi} \log \, 2 \, \sin(\theta_e)
$$
where $L$ is the Lobachevsky function
$$ L(x) =  \int_{0}^{x} \log \, 2 \, \sin \, t \, dt .$$
$N$ is the number of vertices in a fundamental domain and the sum is over the edges in a fundamental domain.
Note that $\sum_{edge \; e} L(\theta_e)$ is the volume of the hyperbolic polyhedron $P$ defined by the dual graph of the isoradial graph(vertices of $P$ are vertices of the dual graph).
And the critical edge weight function and hyperbolic volume do not depend on choices of the common radius of the isoradial graph.
\subsection{Partition function per fundamental domain}
Define the \textit{partition function per fundamental domain} of $G$ to be 
$$\log Z \stackrel{\text{def}}{=} \lim\limits_{n \to \infty} \frac{1}{n^2} \log Z(G_n) .$$
Kenyon, Okounkov and Sheffield proved that\cite[Theorem 3.5]{kenyon2006dimers},
$$ \log Z = \frac{1}{(2\pi i)^2} \int_{\mathbb{T}^2} \log |P(z, w)| \frac{dz}{z} \frac{dw}{w} $$
For isoradial graphs, Tili\`{e}re proved that\cite{tilire2006PartitionFO}
$$
	\log Z = \sum_{edge \; e} \frac{1}{\pi}L(\theta(e))+\frac{\theta(e)}{\pi}\log\,2\,\sin(\theta(e)) 
$$
Combining above equations gives
\begin{equation}\label{eq::vol-mahler formula}
	\begin{split}
		M(P(z, w)) & \stackrel{\text{def}}{=} \frac{1}{(2\pi i)^2} \int_{\mathbb{T}^2} \log |P(z, w)| \frac{dz}{z} \frac{dw}{w} \\
		& = \sum_{edge \; e} \frac{1}{\pi}L(\theta_e)+\frac{\theta_e}{\pi}\log\,2\,\sin(\theta_e) \\
	\end{split}
\end{equation}
for $\mathbb{Z}^2$-periodic isoradial bipartite graphs $G$.

\section{The Champanerkar-Kofman-Lal\'{i}n conjecture} \label{section::CKL conjecture}
In this section we describe the CKL conjecture and prove it for two examples.
\subsection{Construction and statement}\label{section::ckl_construction}
We describe the conjecture following \cite{Champanerkar_2018}.
Let $I = (-1,1)$. 
Let $L$ be a link in the thickened torus $T^2 \times I$ with an alternating diagram on $T^2 \times {0}$, projected onto the 4-valent graph $G_L$. 
The diagram is cellular if the complementary regions are disks, which are called the faces of $L$ or of $G_L$. 
When lifted to the universal cover of $T^2 \times I$, the link $L$ becomes a $ \mathbb{Z}^2 $-periodic alternating link $\mathcal{L}$ in $\mathbb{R}^2 \times I$, such that $L = \mathcal{L}/\Lambda$ for a two-dimensional lattice $\Lambda$ acting by translations of $\mathbb{R}^2$.
The diagram of $L$ on $T^2 \times {0}$ is reduced if four distinct faces meet at every crossing of $G_L$ in $\mathbb{R}^2$.
Link diagrams on $T^2 \times {0}$ will be alternating, reduced and cellular below.

For a $ \mathbb{Z}^2 $-periodic alternating link $\mathcal{L}$ in $\mathbb{R}^2 \times I$, the projection graph $G_{\mathcal{L}}$ in $\mathbb{R}^2$ is $ \mathbb{Z}^2 $-periodic
and can be checkerboard colored. The Tait graph $G_T$ is the planar checkerboard graph for which a vertex is assigned to every shaded region and an edge to every crossing of $\mathcal{L}$.
Then applying the generalized Temperley's trick\cite{kenyonTreesAndMatchings} on $G_T$.
That is taking a barycentric subdivision of $G_T$ and declaring the centers of the edges of $G_T$ to be the white vertices and the centers of the faces of $G_T$ together with the vertices of $G_T$ to be the black vertices.
The resulting graph $G$ is a $ \mathbb{Z}^2 $-periodic bipartite graph.
See Figure \ref{img::triaxialLink} and \ref{img::infWeaveLink}.
Champanerkar-Kofman-Lal\'{i}n chose edge weight 1 for all edges(called \textit{uniform weights} below).
Given edge weights, one can define $P(z, w)$ as in Section \ref{section::background}.

Now we define the bipyramid volume of a link.
Let $B_n$ denote the hyperbolic regular ideal bipyramid whose link polygons at the two coning vertices are regular n-gons.
For a face $f$ of a planar or toroidal graph, let $|f|$ denote the degree of the face; that is, the number of its edges. Let $L$ be an alternating link diagram on the torus as above. Define the bipyramid volume of $L$ as follows:
$$ vol^{\lozenge}(L)=\sum_{f\in \{\text{faces of } L\}} vol(B_{|f|}) $$
where $ vol(B_{|f|}) $ is the hyperbolic volume of  $B_{|f|}$.
Note that $vol(B_2)=0$.
Champanerkar, Kofman and Lal\'{i}n conjectured the following\cite[Conjecture 1]{Champanerkar_2018}.
\begin{conj}
Let $\mathcal{L}$ be any $ \mathbb{Z}^2 $-periodic alternating link, with toroidally alternating quotient link L.
Let $P(z, w)$ be the characteristic polynomial of the toroidal dimer model constructed as above.
Then $$ 2\pi M(P(z, w)) \geq vol^{\lozenge}(L) . $$
\end{conj}

The CKL conjecture and the Vol-Det conjecture are closely related.
The bipyramid volume equals to the hyperbolic volume $vol((T^2 \times I) - L)$ for semi-regular links\cite[Theorem 3.5]{CKP2018}.
Thus the CKL conjecture is equivalent to the toroidal Vol-Det conjecture\cite[Conjecture 6.6]{CKP2018}
$$vol((T^2 \times I) - L) \leq 2\pi M(P(z, w))$$ 
in this case.
And the toroidal Vol-Det conjecture is the asymptotic version of the Vol-Det conjecture in the following sense.
Let $K_n$ be a sequence of alternating links F$\o$lner converging to the $ \mathbb{Z}^2 $-periodic alternating link $\mathcal{L}$, denoted by $K_n \xrightarrow{\text{F}} \mathcal{L}$.
Let $L$ be the alternating quotient link of $\mathcal{L}$ and $P(z, w)$ be the characteristic polynomial of the associated toroidal dimer model.
Champanerkar and Kofman showed that \cite[Theorem 2.4]{champanerkar2016determinant}
$$ K_n \xrightarrow{\text{F}} \mathcal{L} \implies \lim\limits_{n \to \infty} \frac{\log \det (K_n)}{c(K_n)} = \frac{M(P(z, w))}{c(L)} $$
where $ c(\star) $ means the crossing number of $\star$.
Combined with another conjecture of Champanerkar, Kofman and Purcell\cite[Conjecture 6.5]{CKP2018}
\begin{equation} \label{eq::anotherConj}
	\lim\limits_{n \to \infty} \frac{vol(K_n)}{c(K_n)} = \frac{vol((T^2 \times I) - L) }{c(L)} .
\end{equation}
We have that if the Vol-Det conjecture holds for $K_n$, then the toroidal Vol-Det conjecture holds for $L$ assuming Equation \ref{eq::anotherConj}.

\subsection{Proofs for triaxial link and infinite weave link}\label{section::proofForExamples}
\begin{figure}
	\includegraphics[width=\linewidth]{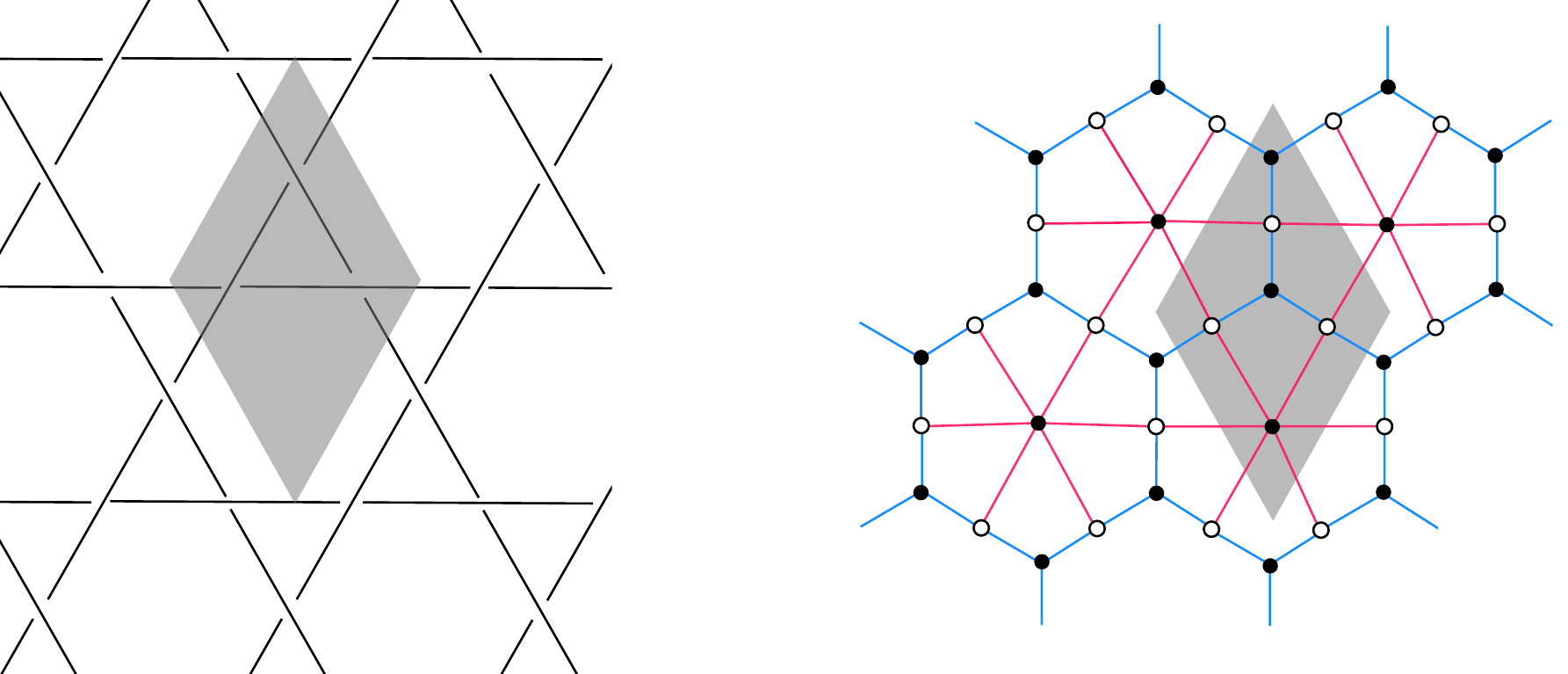}
	\caption{Triaxial link and its induced dimer graph with the fundamental domain.}
	\label{img::triaxialLink}
	
	\includegraphics[width=\linewidth]{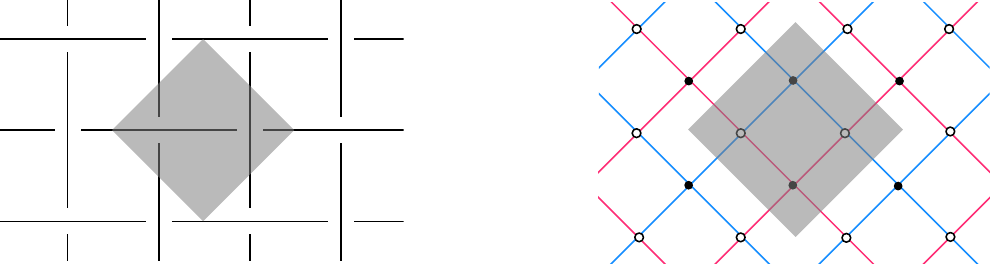}
	\caption{Infinite weave link and its induced dimer graph with the fundamental domain.}
	\label{img::infWeaveLink}
\end{figure}
Champanerkar-Kofman-Lal\'{i}n already proved the conjecture for these two examples using exact computations of Mahler measures\cite[Section 4]{Champanerkar_2018}.
We reprove it using formula \ref{eq::vol-mahler formula}.
Triaxial link and the induced dimer graph(denoted as $\mathcal{D}_T$) are given in Figure \ref{img::triaxialLink} and infinite weave link and the induced dimer graph(denoted as $\mathcal{D}_W$) are given in Figure \ref{img::infWeaveLink}.
The isoradial embeddings of the dimer graphs are given in the figures.

The blue edges of $\mathcal{D}_T$ have rhombus angles $\frac{\pi}{3}$ and weights $\sqrt{3}$.
The red edges of $\mathcal{D}_T$ have rhombus angles $\frac{\pi}{6}$ and weights $1$.
$\mathcal{D}_T$ can be gauged into the uniform edge weights by multiplying all blue edges in a fundamental domain a factor $\frac{1}{\sqrt{3}}$.
Since these edges are exactly the edges incident to two black vertices, the overall effect on the $\log Z$ is $\log \frac{1}{\sqrt{3}} + \log \frac{1}{\sqrt{3}} = -\log 3$. 
Let $P_{T}^{c}(z, w)$ be the partition function of $\mathcal{D}_T$ with critical weights, $P_{T}(z, w)$ be that of $\mathcal{D}_T$ with uniform weights.
Then $ M(P_{T}(z, w)) = M(P_{T}^{c}(z, w)) - \log\,3 $.
By formula \ref{eq::vol-mahler formula},
\begin{equation}
	\begin{split}
		M(P_{T}^{c}(z, w))   & = \frac{6}{\pi} \Lambda(\frac{\pi}{3})+\frac{6}{\pi} \Lambda(\frac{\pi}{6}) +6\times\frac{1}{3}\log\sqrt{3}  \\
		& = \frac{6}{\pi} \times \frac{1}{3}v_{tet} + \frac{6}{\pi} \times \frac{1}{2}v_{tet}+ \log\,3 \\
		& = \frac{5}{\pi} v_{tet}+ \log\,3 
	\end{split}
\end{equation}
where $v_{tet}$ is the volume of the regular ideal hyperbolic tetrahedron.
Thus $ 2\pi M(P_{T}(z, w)) = 10 v_{tet} $.
Since the bipyramid volume of the triaxial link is $2 vol(B_3)+vol(B_6) = 10 v_{tet}$, the above calculation proves that the CKL conjecture is true for the triaxial link and the equality is obtained.

The edges of $\mathcal{D}_W$ has rhombus angles $\frac{\pi}{2}$ and weights $\sqrt{2}$.
$\mathcal{D}_W$ can be gauged into a dimer model with uniform weights by multiplying all edges in a fundamental domain a factor $\frac{1}{\sqrt{2}}$.
Since these edges are exactly the edges incident to two black vertices, the overall effect on the $\log Z$ is $\log \frac{1}{\sqrt{2}} + \log \frac{1}{\sqrt{2}} = -\log 2$. 
Let $P_{W}^{c}(z, w)$ be the partition function of $\mathcal{D}_W$ with critical weights, $P_{W}(z, w)$ be that of $\mathcal{D}_W$ with uniform weights.
Then $ M(P_{W}(z, w)) = M(P_{W}^{c}(z, w)) - \log\,2 $.
By formula \ref{eq::vol-mahler formula},
\begin{equation}
	\begin{split}
	 M(P_{W}^{c}(z,w)) & = \frac{8}{\pi}\Lambda(\frac{\pi}{4})+8\times\frac{1}{4}\log\sqrt{2}  \\
		& = \frac{8}{\pi}\times(\frac{1}{8}v_{oct})+ \log\,2 \\
		& = \frac{1}{\pi} v_{oct} + \log\,2 
	\end{split}
\end{equation}
where $v_{oct}$ is the volume of the regular ideal hyperbolic octahedron.
Thus $ 2\pi M(P_{W}(z, w)) = 2 v_{oct} $.
Since the bipyramid volume of the infinite weave link is $2 v_{oct}$, the above calculation proves that the CKL conjecture is true for the infinite weave link and the equality is obtained.

These two links have special properties, see the work of Champanerkar, Kofman and Purcell\cite{CKP2018} and Kaplan-Kelly\cite{Kaplan-Kelly}.

\section{Spanning tree entropy}\label{section::spanning tree entropy}
In this section, we verify the conjecture using the work of Champanerkar and Kofman\cite{champanerkar2016determinant} and some calculations of Teufl and Wagner\cite{teufl2010number} for some examples.

Let $G$ be the associated toroidal dimer model of the $ \mathbb{Z}^2 $-periodic link $\mathcal{L}$ and $T$ be one of the Tait graph of $\mathcal{L}$.
By a result of Kenyon, Okounkov and Sheffield\cite{kenyon2006dimers} 
$$ M(P(z, w)) = \log Z(G) .$$
And Champanerkar and Kofman showed that\cite[Proposition 5.3]{champanerkar2016determinant} 
$$ \log Z(G) = \lim\limits_{n \to \infty} \frac{\log Z(G_n)}{n^2} = \lim\limits_{n \to \infty} \frac{\log N_{ST} (T_{n})}{n^2} $$
where $G_{n} = G/(n\Lambda)$ is a finite quotient of $G$, $T_n = T \cap (n\Lambda)$ and $N_{ST}(T)$ denotes the number of spanning trees of $T$.

Denote the spanning tree entropy of the lattice $T$ by $z_{T}$
$$ z_{T} =  \lim\limits_{n \to \infty} \frac{\log N_{ST} ({T}_{n})}{V(T_n)} $$
where $V(\star)$ is the number of vertices of $\star$.
Denote the spanning tree entropy per fundamental domain of the lattice $T$ by $z^{fd}_{T}$
$$ z^{fd}_{T} = \lim\limits_{n \to \infty} \frac{\log N_{ST} (T_{n})}{n^2} .$$
We have
$$  z^{fd}_{T} = N_v \, z_{G_{T}} . $$
where $N_v$ is the number of vertices of $G_T$ in the fundamental domain.
Combining these equations, we have
$$ M(P_{G}(z, w)) = z^{fd}_{T}  .$$

Now we prove the CKL conjecture for some examples using this formula.
Example 1 and Example 4 are already calculated by Champanerkar, Kofman and Lat\'{i}n\cite{Champanerkar_2018}.
The approximate values of volume of bipyramids comes from Adams\cite[Table 1]{adamsbipyramids}.
$P$ denotes the characteristic polynomial of the associated toroidal dimer model of the $ \mathbb{Z}^2 $-periodic link $\mathcal{L}$.

\begin{remark} \label{remark::error}
	Teufl and Wagner\cite[Section 3]{teufl2010number} calculated the spanning tree entropy per vertice of a lattice in terms of another.
	However, the computation does not directly apply to spanning tree entropy per fundamental domain which we need to calculate the Mahler measure.
	We adjust their formulas for spanning tree entropy per fundamental domain in the following.
\end{remark}

\begin{figure}
	\begin{subfigure}{0.5\textwidth}
		\centering
		\includegraphics[width=0.7\linewidth]{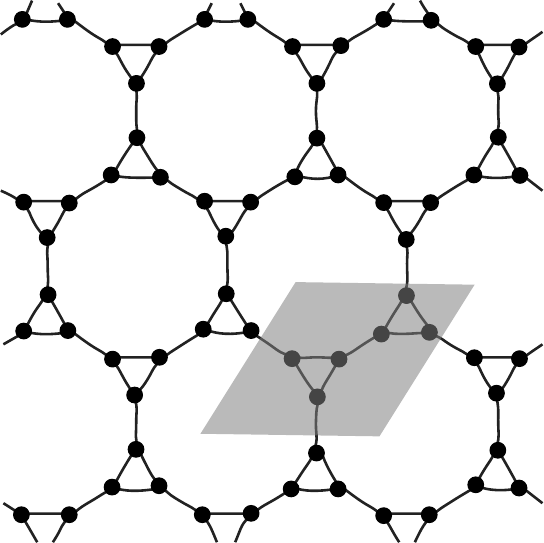}
		\caption{lattice $T_{3 \cdot 12 \cdot 12}$}
		\label{img::triak}
	\end{subfigure}%
	\begin{subfigure}{0.5\textwidth}
		\centering
		\includegraphics[width=0.7\linewidth]{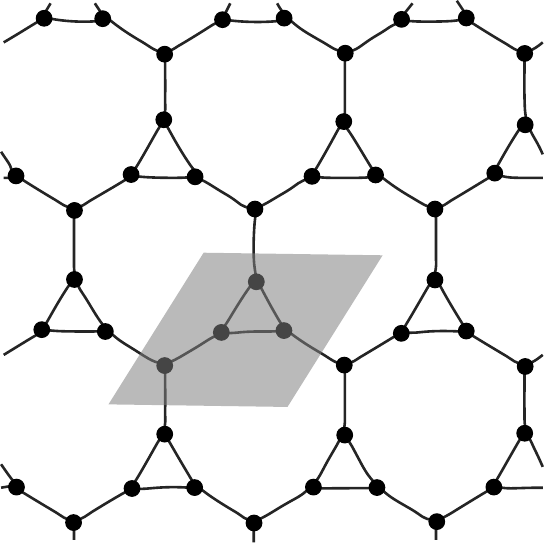}
		\caption{lattice $T_{nine}$}
		\label{img::nine}
	\end{subfigure}
	\begin{subfigure}{0.5\textwidth}
		\centering
		\includegraphics[width=0.7\linewidth]{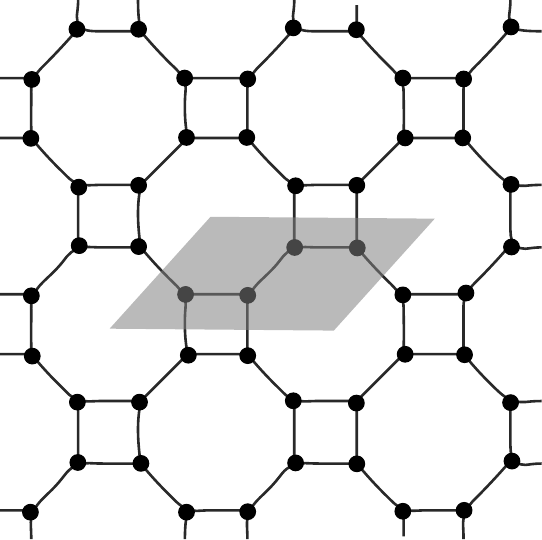}
		\caption{lattice $T_{4 \cdot 8 \cdot 8}$}
		\label{img::bathroom}
	\end{subfigure}
	\begin{subfigure}{0.5\textwidth}
		\centering
		\includegraphics[width=0.7\linewidth]{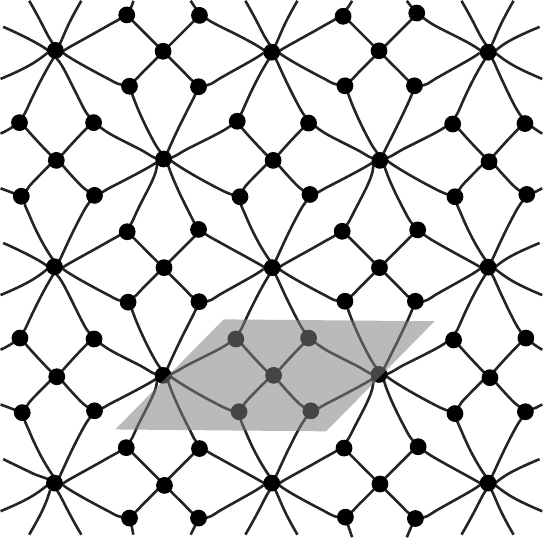}
		\caption{lattice $T_{kite}$}
		\label{img::kite}
	\end{subfigure}
	\caption{}
\end{figure}

\begin{description}
	\item[Example 1.] 
	One of the Tait graphs of the $ \mathbb{Z}^2 $-periodic Rhombitrihexagonal link\cite[Section 4.3]{Champanerkar_2018} is the kagom\'{e} lattice $T_{kag}$ in \cite{teufl2010number}.
	Its spanning tree entropy per fundamental domain is
	$$ z^{fd}_{kag} = z^{fd}_{tri} + \log 6 \approx 3.407088 .$$
	where $ z^{fd}_{tri} = z_{tri} \approx 1.615329 $ is the spanning tree entropy per fundamental domain of the triangular lattice\cite{teufl2010number}.
	Thus $$ 2\pi M(P_{kag}) = 2\pi \times z^{fd}_{kag} \approx 21.407368  $$
	The result agrees with the calculations in the proof of Corollary 14 in \cite{Champanerkar_2018}.
	
	\item[Example 2.] 
	See Figure \ref{img::triak}.
	The spanning tree entropy per fundamental domain of $3 \cdot 12 \cdot 12$ lattice is\cite{teufl2010number}
	$$ z^{fd}_{3 \cdot 12 \cdot 12} = z^{fd}_{tri} + \log 15 \approx 4.323379 .$$
	Thus $$ 2\pi M(P_{3 \cdot 12 \cdot 12}) = 2\pi \times  z^{fd}_{3 \cdot 12 \cdot 12} \approx 27.164592  $$
	And
	\begin{equation}
		\begin{split}
			vol^{\lozenge}(L_{3 \cdot 12 \cdot 12}) & =  vol(B_{12}) + 8 \, vol(B_3)   \\
			&\approx 10.3725 + 8 \times 2.0298 \\
			& = 26.6109 \\
			& < 2\pi M(P_{3 \cdot 12 \cdot 12}) . \\
		\end{split}
	\end{equation}
	
	\item[Example 3.] 
	See Figure \ref{img::nine}.
	By \cite{teufl2010number}
	$$ z^{fd}_{nine} =  z^{fd}_{tri} + 2 \log 2 \approx 3.001623 .$$
	Thus $$ 2\pi M(P_{nine}) = 2\pi \times z^{fd}_{nine} \approx 18.859756  $$
	And
	\begin{equation}
		\begin{split}
			vol^{\lozenge}(L_{nine}) & =  vol(B_9) + 5 \, vol(B_3)   \\
			&\approx 8.5836 + 5 \times 2.0298 \\
			& = 18.7326 \\
			& < 2\pi M(P_{nine}) . \\
		\end{split}
	\end{equation}
	
	\item[Example 4.] 
	See Figure \ref{img::bathroom}.
	Chang and Shrock gave\cite[Section 5]{chang2006some} a closed formula for the spanning tree entropy of $4 \cdot 8 \cdot 8$ lattice $T_{4 \cdot 8 \cdot 8}$ and gave a approximate value with high accuracy:
	$$ z_{4 \cdot 8 \cdot 8} \approx 0.786684275378832 .$$
	There are 4 vertices in the fundamental domain of $T_{4 \cdot 8 \cdot 8}$.
	Thus $$ 2\pi M(P_{4 \cdot 8 \cdot 8}) = 2\pi \times 4 \times z_{4 \cdot 8 \cdot 8} \approx 19.7715323218  $$
	And
	\begin{equation}
		\begin{split}
			vol^{\lozenge}(L_{4 \cdot 8 \cdot 8}) & = vol(B_8) +  vol(B_4) + 4 \, vol(B_3)  \\
										& \approx 7.8549 + 3.6638 + 4 \times 2.0298 \\
										& = 19.6379 \\
										& < 2\pi M(P_{4 \cdot 8 \cdot 8}).  \\
		\end{split}
	\end{equation}
	The result agrees with the calculations in the proof of Theorem 19 in \cite{Champanerkar_2018}.
	See the Appendix for the relation of two closed formulas.
	
	\item[Example 5.] 		
	See Figure \ref{img::kite}.
	By \cite{teufl2010number}
	$$ z^{fd}_{kite} = z^{fd}_{4 \cdot 8 \cdot 8} + 2\log 6 \approx 6.730256 .$$
	Thus $$ 2\pi M(P_{kite}) = 2\pi \times z^{fd}_{kite} \approx 42.287446  $$
	And
	\begin{equation}
		\begin{split}
			vol^{\lozenge}(L_{kite}) & = vol(B_8) +  7 \, vol(B_4) + 4 \, vol(B_3)  \\
			& \approx 7.8549 + 7 \times 3.6638 + 4 \times 2.0298 \\
			& = 41.6207  \\
			& < 2\pi M(P_{kite}).  \\
		\end{split}
	\end{equation}

\end{description}

\appendix
\section{}
In this appendix, we show the closed formula of $ 2\pi M(P_{4 \cdot 8 \cdot 8}) $ given by Champanerkar, Kofman and Lat\'{i}n and the closed formula of $2\pi \times 4 \times z_{4 \cdot 8 \cdot 8} $ given by Chang and Shrock are equal.

Chang and Shrock showed that
$$ z_{4 \cdot 8 \cdot 8} = \frac{C}{\pi} + \frac{1}{2} \log (\sqrt{2}-1) + \frac{1}{\pi} \text{Ti}_2 (3+2\sqrt{2}) $$
where $C= D(i) = vol_{oct}/4$ is the Catalan number and $\text{Ti}_2$ is the inverse tangent function.
Thus
\begin{equation}
	\begin{split}
		2\pi \times 4 \times z_{4 \cdot 8 \cdot 8} 
		& = 2\pi \times 4 \times (\frac{C}{\pi} + \frac{1}{2} \log (\sqrt{2}-1) + \frac{1}{\pi} \text{Ti}_2 (3+2\sqrt{2})) \\ 
		& = 8C + 4\pi \log (\sqrt{2}-1) + 8 \text{Ti}_2 (3+2\sqrt{2}) \\ 
	\end{split}
\end{equation}

Champanerkar, Kofman and Lat\'{i}n showed that 
\begin{multline}
	$$ 2\pi M(P_{4 \cdot 8 \cdot 8}) 
	= 8D(i) + \arccos (-\frac{7}{9}) \log (17+12\sqrt{2}) \\ 
	+ 4D(\frac{\sqrt{7+4\sqrt{2}i}}{3})-4D(-\frac{\sqrt{7+4\sqrt{2}i}}{3}) $$
\end{multline}
where $$ D(z) = \text{Im}(\text{Li}_2 (z)) + \arg(1-z) \log|z| $$ is the Bloch-Wigner dilogarithm and arg denotes the branch of the argument lying between $-\pi$ and $\pi$\cite[Section 3]{zagier_dilogarithm}.

\begin{figure}
	\begin{subfigure}{0.33\textwidth}
		\centering
		\includegraphics[width=1.1\linewidth]{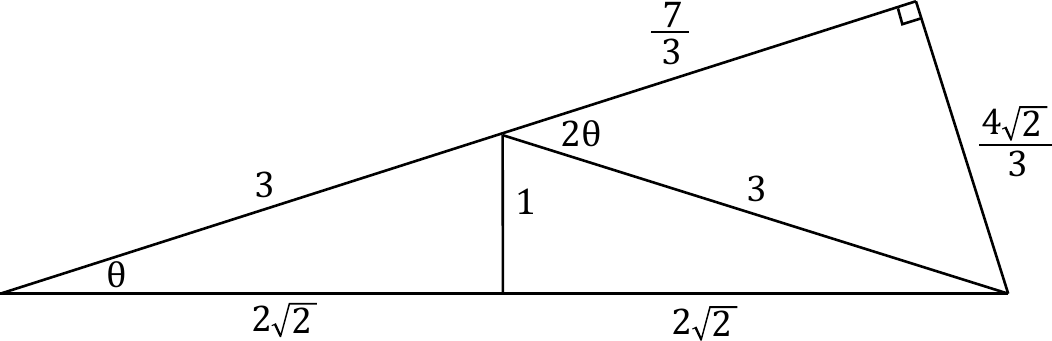}
		\caption{}
		\label{img::triangle1}
	\end{subfigure}%
	\begin{subfigure}{0.33\textwidth}
		\centering
		\includegraphics[width=0.4\linewidth]{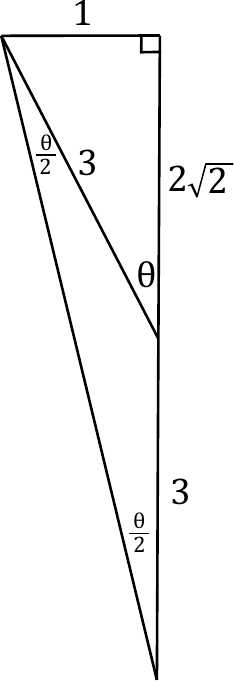}
		\caption{}
		\label{img::triangle2}
	\end{subfigure}%
	\begin{subfigure}{0.33\textwidth}
		\centering
		\includegraphics[width=1.1\linewidth]{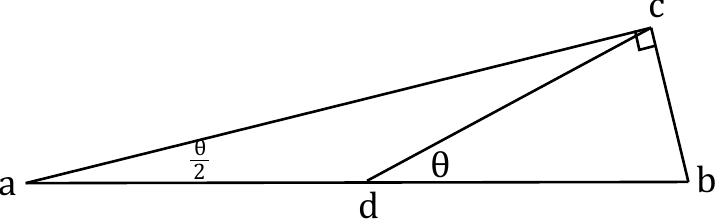}
		\caption{}
		\label{img::triangle3}
	\end{subfigure}
	\caption{}
\end{figure}

Here 
$$ \log (17+12\sqrt{2})  = 2 \log (3+2\sqrt{2})  = 4 \log(1+\sqrt{2})  $$
and 
$$ \frac{\sqrt{7+4\sqrt{2}i}}{3} = \frac{2\sqrt{2}+i}{3} = e^{i\theta} $$
$$ \arccos (-\frac{7}{9}) = \pi - 2\theta $$
by Figure \ref{img::triangle1}.

By \cite[Chapter 2]{lewinPolylogarithms}, the inverse tangent integral $\text{Ti}_2 (x) = \int_{0}^{x} \frac{tan^{-1}(t)}{t} dt$ satisfies that
$$ \text{Ti}_2 (x) = \text{Im}(\text{Li}_2 (ix)) \quad \quad  x\in \mathbb{R} $$
where 
$$ \text{Li}_2(z) = - \int_{0}^{z} \log(1-u) \frac{du}{u} \quad \quad   z \in \mathbb{C} \backslash [1, \infty) $$ 
is the analytic continuation of the dilogarithm function\cite{zagier_dilogarithm}.
Thus 
$$ \text{Ti}_2 (3+2\sqrt{2}) = D((3+2\sqrt{2})i) - \arg(1-(3+2\sqrt{2})i) \log |(3+2\sqrt{2})i| .$$
Here 
$$ \arg(1-(3+2\sqrt{2})i) = -(\frac{\pi}{2} - \frac{\theta}{2}) $$
by Figure \ref{img::triangle2}.
So
\begin{equation}
	\begin{split}
		 \arg(1-(3+2\sqrt{2})i) \log |(3+2\sqrt{2})i| 
		& =  -(\frac{\pi}{2} - \frac{\theta}{2})  \log |(3+2\sqrt{2})i| \\
		& = -(\pi - \theta) \log(1+\sqrt{2}) \\ 
	\end{split}
\end{equation}
and
\begin{equation}
	\begin{split}
		4\pi & \log (\sqrt{2}-1)  + 8 \text{Ti}_2 (3+2\sqrt{2}) \\
		& = 4\pi \log (\sqrt{2}-1) + 8D((3+2\sqrt{2})i) + 8(\pi - \theta) \log(1+\sqrt{2})\\
		& = 8D((3+2\sqrt{2})i) + 4(\pi-2\theta)\log(1+\sqrt{2}) \\
		& = 8D((3+2\sqrt{2})i) + \arccos (-\frac{7}{9}) \log (17+12\sqrt{2}) \\
	\end{split}
\end{equation}
It remains to show that 
\begin{equation} \label{eq::volumeIdentity}
 8D((3+2\sqrt{2})i) = 4D(e^{i\theta})-4D(-e^{i\theta})  
\end{equation}
in order to prove their formulas are equal.
$D(z)$ is the volume of the ideal hyperbolic tetrahedron with vertices $0$, $1$, $z$ and $\infty$ where $z$ is in the upper half plane.
When $z$ is in the lower half plane, $D(z) = -D(\overline{z})$.
Using the identity $ D(z) = -D(\frac{1}{z}) $, we have
$$ 4D(e^{i\theta})-4D(-e^{i\theta}) = 4D(e^{i\theta}) + 4D(e^{i(\pi-\theta)}) . $$
$D(e^{i\theta})$ is the volume of the ideal tetrahedron $\Delta_1$ with vertices $d, b, c, \infty$ in Figure \ref{img::triangle3}.
$D(e^{i(\pi-\theta)})$ is the volume of the ideal tetrahedron $\Delta_2$ with vertices $d, c, a, \infty$ in Figure \ref{img::triangle3}.
$D((3+2\sqrt{2})i)$ is the volume of the ideal tetrahedron $\Delta_3$ with vertices $a, b, c, \infty$ by Figure \ref{img::triangle2} and \ref{img::triangle3}.
Now by the reflection isometry along the hemisphere defined by the vertices $a, b, c$ in Figure \ref{img::triangle3} of the tetrahedra $\Delta_1$ and $\Delta_2$, we have Equation \ref{eq::volumeIdentity}.

\bibliographystyle{plain}
\bibliography{References}

\end{document}